\pdfoutput=1
\documentclass[12pt, reqno]{amsart}
\usepackage{amsmath}
\usepackage[T1,T4]{fontenc}
\usepackage{textgreek}
\usepackage{mathrsfs}
\allowdisplaybreaks
\begin{document}
\setcounter{page}{1}

\title[\hfilneg \hfil Fractional Integro-differential Equations]
{On Implicit $\pmb{\psi}$-Caputo Fractional Integro-differential Equations}

\author[Deepak Pachpatte \hfil]
{Deepak B. Pachpatte}

\address{Deepak B. Pachpatte \newline
 Dept. of Mathematics,
 Dr. B. A. M. University, Aurangabad,
 Maharashtra 431004, India}
\email{pachpatte@gmail.com}

\author[Bhagwat Yewale \hfil]
{Bhagwat R. Yewale}

\address{Bhagwat R. Yewale \newline
 Dept. of Mathematics,
 Dr. B. A. M. University, Aurangabad,
 Maharashtra 431004, India}
\email{yewale.bhagwat@gmail.com}


\begin{abstract}
  In this paper, we investigate the existence and uniqueness of solution of nonlinear
  $\psi$-Caputo fractional differential equation with the help of Banach fixed point theorem.
  Moreover, by using $\psi$-Gronwall inequality, we studied some qualitative properties of
  solutions such as  estimate and stability.
\end{abstract}

\maketitle

\section{Introduction}
Fractional calculus is widely acknowledged as one of the important fields behind the development
of theory of integration and differentiation as it generalizes the integer order integral and differential
to the non-integer order. The study of fractional calculus occupies a very privileged position in the field
of modern research due to its applicability in modeling many physical phenomena, for example,
Medical and Health Sciences \cite{Kum}, Image edge detection \cite{Nan}, Bioengineering \cite{Mag}.
In fractional calculus, there has been considerable interest in developing various
fractional integral and differential operators. Starting from the definition of Riemann-Liouville
fractional integral and differential operators several new fractional integrals and differential operators
have been introduced  such as Caputo, Hadamard, Hilfer and many more can be found in the
literature \cite{Kil}. Recently, with the development of the theory of fractional calculus, many authors proposed
fractional integrals and derivatives of a function with respect to another function known as $\psi$-Fractional Calculus.
For instance, Almeida \cite{Alm1}, introduced a fractional derivative of a function with respect to another
function known as $\psi$-Caputo fractional derivative. Based on Hilfer fractional derivative, Sousa and Oliveira \cite{Sou},
established a new fractional derivative called  $\psi$-Hilfer fractional derivative. In \cite{Alm, Alm2, Pac1, Sam, Wah}, authors have
studied various inequalities, fractional differential equations, and their properties using $\psi$-Fractional operators.

\par In recent years, authors devoted themselves to study the existence and uniqueness of fractional differential equations.
For instance, Abdo et al. \cite{Abd}, studied the existence and uniqueness of the solution of $\psi$-Hilfer fractional
integro-differential equations by employing Banach and Krasnoselskii's fixed point theorem. Sousa and Kucche \cite{Sou1},
investigated the existence, uniqueness, and stability of solutions of $\psi$-Caputo impulsive fractional functional
differential equations through fixed point theorem. Pachpatte\cite{Pac2}, studied the existence, uniqueness and properties of
Fredholm-type fractional integro-differential for $\psi$-Hilfer fractional derivative. With the growing range of research, the theory of fractional calculus progressed rapidly and acknowledged as a significant tool to describe qualitative analysis of solutions of
fractional differential equations such as existence, uniqueness, estimate, stability, continuous dependence.
For more details see \cite{Der, Fur, Kuc1, Kuc2,Pac3, Tid}.

Motivated by the above work, in this paper, we consider the ${\psi}$-Caputo fractional integro-differential equation
of the type
\[\label{eq1.1}
\textrm{\m{N}}^{\alpha, \psi}_{0+}\vartheta(\textrm{\m{z}})
=\mathscr{F}
  \big(\textrm{\m{z}}, \vartheta(\textrm{\m{z}}),
    \textrm{\m{I}}^{\alpha, \psi}_{0+}
    \mathscr{H}
       (\textrm{\m{z}}, \tau, \textrm{\m{N}}^{\alpha, \psi}_{0+}\vartheta(\tau))
  \big),
       \;\mathrm{for\; all}\; \textrm{\m{z}}\in [0, b], 0<\alpha<1,
\tag{1.1}
\]
\[\label{eq1.2}
\vartheta(0)=\vartheta_{0}, \qquad \vartheta_{0}\in\mathbb{R}
\tag{1.2}
\]
where $\textrm{\m{N}}^{\alpha, \psi}_{0+}$ denotes the $\psi$-Caputo fractional derivative and
$\textrm{\m{I}}^{\alpha, \psi}_{0+}$ denotes the $\psi$-Riemann-Liuoville fractional integral.
$\mathscr{F}:[0, b]\times\mathbb{R}\times\mathbb{R}\rightarrow\mathbb{R}$ and
$\mathscr{H}:[0, b]\times[0,b]\times\mathbb{R}\rightarrow\mathbb{R}$ are continuous functions satisfying certain assumptions.

 This paper is organized as follows: In section 2, we present some preliminaries and prove some basic lemmas which
 are used to carry out our work. In section 3, we prove existence and uniqueness of solution of the
 problem \eqref{eq1.1}-\eqref{eq1.2}. Estimate on solutions and  different aspects of stability results
 of solutions of \eqref{eq1.1}-\eqref{eq1.2} are derived in section 4 and section 5 respectively.
 Finally, illustrative example is given to demonstrate our results.

\section{Preliminaries}
In this section, we give some basic definitions, preliminaries which are useful for our subsequent discussions.
Let $\textbf{J}=[0, b] \subseteq \mathbb{R^{+}}$.
Let $\mathcal{C}(\textbf{J},\mathbb{R})$ be the Banach space of all continuous functions $\vartheta:\textbf{J}\rightarrow \mathbb{R}$
with
$\|\vartheta\|= \underset{\textrm{\m{z}}\in\textbf{J}}{\mathrm{sup}}{|\vartheta(\textrm{\m{z}})|},$
$\mathcal{L}_{p}(\textbf{J})$, $1\leq p \leq \infty$ be the space of Lebesgue measurable functions endowed with the norm
$\|\vartheta\|_{\mathcal{L}_{p}(\textbf{J})}<\infty$
and $\mathcal{C}^{n}(\textbf{J})$ be the space of n-times continuously differentiable functions on \textbf{J}.

Now we highlight some definitions of a fractional integrals and derivative of a function with respect to another function:\newline
{\textbf{Definition 2.1}}. \cite{Kil} Let $\alpha>0,$ $\vartheta\in \mathcal{L}_{1}(\textbf{J})$ and ${\psi}$ be an increasing and positive monotone function on \textbf{J}, having a continuous derivative $\psi^{'}$ on \textbf{J}. Then $\psi$-Riemann Liouville fractional integral of $\vartheta$ with respect to $\psi$ of order $\alpha$ is given by
\[
\textrm{\m{I}}^{\alpha, \psi}_{0+}\vartheta(\textrm{\m{z}})
=\frac{1}{\Gamma(\alpha)}
\int_{0}^{\textrm{\m{z}}}{\psi}^{'}(\textkappa)({\psi}(\textrm{\m{z}})-{\psi}(\textkappa)^{\alpha-1}
\vartheta(\textkappa)d\textkappa,\;\textrm{\m{z}}>0.
\]
{\textbf{Definition 2.2}}. \cite{Alm1} Let $\alpha>0$ and $\vartheta$, $\psi \in\mathcal{C}^{n}(\textbf{J})$ with $\psi$ an increasing function such that $\psi^{'}(\textrm{\m{z}})\neq 0,$ for all $\textrm{\m{z}} \in \textbf{J}$, the $\psi$-Caputo fractional derivative of a function $\vartheta$ with respect to ${\psi}$ of order $\alpha$ is defined as
\[
\textrm{\m{N}}^{\alpha,\psi}_{0+}\vartheta(\textrm{\m{z}})
=\frac{1}{\Gamma(n-\alpha)}\int_{0}^{\textrm{\m{z}}}\psi^{'}(\textkappa)(\psi(\textrm{\m{z}})-\psi(\textkappa))^{n-\alpha-1}
\bigg(\frac{1}{\psi^{'}(\textrm{\m{z}})}\frac{d}{d\textrm{\m{z}}}\bigg)^{n}\vartheta(\textkappa)d\textkappa,
\]
where $n=\lceil \alpha \rceil+1$, $\lceil \alpha \rceil$ is the integer part of $\alpha$.\newline

{\textbf{Lemma 2.1}}. \cite{Alm1} Let $n-1< \alpha <n$ and $\vartheta \in \mathcal{C}^{n}(\textbf{J})$ then
\[
\textrm{\m{I}}^{\alpha, \psi}_{0+}\textrm{\m{N}}^{\alpha,\psi}_{0+}\vartheta(\textrm{\m{z}})
=\vartheta(\textrm{\m{z}})-\sum_{m=0}^{n-1}\frac{\vartheta^{[m]}_{\psi}(0)}{m!}(\psi(\textrm{\m{z}})-\psi(0))^{m},
\]
where $\vartheta^{[m]}_{\psi}(\kappa)=\big(\frac{1}{\psi^{'}(\kappa)} \frac{d}{d\kappa}\big)^{m}\vartheta(\kappa)$.

{\textbf{Lemma 2.2}}. \cite{Alm1}
Let $\alpha>0$ and $\vartheta \in \mathcal{C}^{1}(\textbf{J},\mathbb{R}),$ we have
\[
\textrm{\m{N}}^{\alpha,\psi}_{0+}\textrm{\m{I}}^{\alpha, \psi}_{0+}\vartheta(\textrm{\m{z}})
=\vartheta(\textrm{\m{z}}).
\]

{\textbf{Lemma 2.3}}. \cite{Kil} Let $\alpha,\; \beta>0$ and  $\vartheta\in \mathcal{L}_{1}(\textbf{J})$, we have
\[
\textrm{\m{I}}^{\alpha, \psi}_{0+}\textrm{\m{I}}^{\beta, \psi}_{0+}\vartheta(\textrm{\m{z}})
=\textrm{\m{I}}^{\alpha+\beta, \psi}_{0+}\vartheta(\textrm{\m{z}}), \quad \textrm{\m{z}}>0.
\]

{\textbf{Lemma 2.4}}. \cite{Gra} Let $(\Delta, d)$ be complete metric space and $\eta:\Delta\rightarrow\Delta$ be a contracting mapping.
Then $\eta$ has unique fixed point in ${\Delta}.$

{\textbf{Lemma 2.5}} \cite{Sou2}
Let $\vartheta, \eta$ be two integrable functions and $\rho$ be continuous function on \textbf{J}.
Let $\psi\in \mathcal{C}(\textbf{J},\mathbb{R})$ be an increasing function with $\psi^{'}(\textrm{\m{z}})\neq0,$ for $\textrm{\m{z}}\in \textbf{J}$. Assume that\\
(i)$\vartheta$ and $\eta$ are nonnegative,\\
(ii)$\rho$ is nonnegative and nondecreasing.\\
If
\[
\vartheta(\textrm{\m{z}})\leq \eta(\textrm{\m{z}})+\rho(\textrm{\m{z}})\int_{0}^{\textrm{\m{z}}}\psi^{'}(\textkappa)
(\psi(\textrm{\m{z}})-\psi(\kappa))^{\alpha-1}\vartheta(\kappa)d\kappa,
\]
then
\[
\vartheta(\textrm{\m{z}})\leq \eta(\textrm{\m{z}})\int_{0}^{\textrm{\m{z}}}
\sum_{m=1}^{\infty}\frac{[\rho(\textrm{\m{z}})\Gamma(\alpha)]^{m}}{\Gamma(\alpha m)}
\psi^{'}(\kappa)
(\psi(\textrm{\m{z}})-\psi(\kappa))^{\alpha-1}\eta(\kappa)d\kappa,
\]
for $\textrm{\m{z}}\in\textbf{J}$.

{\textbf{Remark 2.1}}\cite{Sou2}
Under the assumptions of Lemma (2.5), let $\vartheta(\textrm{\m{z}})$ be a nondecreasing function on
\textbf{J}. Then we have
\[
\vartheta(\textrm{\m{z}})\leq \eta(\textrm{\m{z}})E_{\alpha}(\rho(\textrm{\m{z}})
\Gamma(\alpha)(\psi(\textrm{\m{z}})-\psi(0))^{\alpha}),
\]
where $E_{\alpha}(\textrm{\m{z}})=\sum_{m=0}^{\infty}\frac{\textrm{\m{z}}^{\alpha}}{\Gamma(\alpha+1)}.$

{\textbf{Lemma 2.6}}
Let $0<\alpha<1$ and $\mathscr{F}:\textbf{J}\times \mathbb{R}\times \mathbb{R}\rightarrow \mathbb{R}$ be a continuous function.
Then the problem \eqref{eq1.1}-\eqref{eq1.2} is equivalent to
\[\label{eq2.1}
\vartheta(\textrm{\m{z}})
=\vartheta_{0}
+\textrm{\m{I}}^{\alpha, \psi}_{0+}
\mathscr{F}
  \big(\textrm{\m{z}}, \vartheta(\textrm{\m{z}}),
    \textrm{\m{I}}^{\alpha, \psi}_{0+}
    \mathscr{H}
       (\textrm{\m{z}}, \tau, \textrm{\m{N}}^{\alpha, \psi}_{0+}\vartheta(\tau))
   \big).
\tag{2.1}\]
{\textbf{Proof.}}
Assume that $\vartheta(\textrm{\m{z}})$ is the solution of \eqref{eq1.1}-\eqref{eq1.2}.
Operating $\textrm{\m{I}}^{\alpha, \psi}_{0+}$ on both the sides of \eqref{eq1.1} and then by using Lemma (2.1), we obtain

\[
 \vartheta(\textrm{\m{z}})-\vartheta(0)=\textrm{\m{I}}^{\alpha, \psi}_{0+}
\mathscr{F} \big(\textrm{\m{z}}, \vartheta(\textrm{\m{z}}),
    \textrm{\m{I}}^{\alpha, \psi}_{0+}
    \mathscr{H}
       (\textrm{\m{z}}, \tau, \textrm{\m{N}}^{\alpha, \psi}_{0+}\vartheta(\tau))\big)
\]
By using \eqref{eq1.2}, we get \eqref{eq2.1}.\newline
Conversely, suppose that $\vartheta$ satisfies the fractional integral equation \eqref{eq2.1}.
Then clearly $\vartheta(0)=\vartheta_{0}$ and using lemma (2.2), we get \eqref{eq1.1}.$\qed$

\section{Existence and Uniqueness}
This section deals with the existence and uniqueness of solution of the problem.\newline
{\textbf{Theorem 3.1}}
Assume that the following hypothesis hold:

($\mathscr{Q1}$) The function $\mathscr{F}$ in \eqref{eq1.1} satisfies following condition
\[
| \mathscr{F}(\textrm{\m{z}}, \vartheta_{1},\vartheta_{2})
-\mathscr{F}(\textrm{\m{z}}, \bar{\vartheta_{1}},\bar{\vartheta_{2}})|
\leq \textrm{\m{n}}(\textrm{\m{z}})
\big[|\vartheta_{1}-\bar{\vartheta_{1}}|+| \vartheta_{2}-\bar{\vartheta_{2}}|\big],
\]
where $\textrm{\m{n}}(\textrm{\m{z}})\in \mathcal{C}(\textbf{J},\mathbb{R})$.

($\mathscr{Q2}$) The function $\mathscr{H}$ in \eqref{eq1.1} satisfies following condition
\[
|\mathscr{H}(\textrm{\m{z}}, \tau,\vartheta)
- \mathscr{H}(\textrm{\m{z}}, \tau, \bar{\vartheta})|
\leq \textrm{\m{k}}(\textrm{\m{z}})
|\vartheta-\bar{\vartheta}|,
\]
where $\textrm{\m{k}}(\textrm{\m{z}})\in \mathcal{C}(\textbf{J},\mathbb{R})$.\\
Define
\begin{align*}
 \mathscr{W}_{\textbf{1}}= \textrm{sup}\{|\textrm{\m{n}}(\textrm{\m{z}})|: 0<\textrm{\m{z}}<b\},\\
\mathscr{W}_{\textbf{2}}= \textrm{sup}\{|\textrm{\m{k}}(\textrm{\m{z}})|: 0<\textrm{\m{z}}<b\}. \\
\end{align*}
If
\[\label{eq3.1}
\frac{(\psi(b)-\psi(0))^{\alpha}}{\Gamma(\alpha+1)}
\mathscr{W}_{\textbf{1}}(1+\mathscr{W}_{\textbf{1}})<1,
\tag{3.1}\]
then there exists a unique solution for the problem \eqref{eq1.1}-\eqref{eq1.2}.\\
{\textbf{Proof.}}
Define \[
(\textrm{\m{P}}\vartheta)(\textrm{\m{z}})
=\vartheta_{0}+\textrm{\m{I}}^{\alpha, \psi}_{0+}
\mathscr{F}\big(\textrm{\m{z}}, \vartheta(\textrm{\m{z}}), \textrm{\m{I}}^{\alpha, \psi}_{0+}
    \mathscr{H}(\textrm{\m{z}}, \tau, \textrm{\m{N}}^{\alpha, \psi}_{0+}\vartheta(\tau))\big).
\]
Then for $\vartheta, \vartheta^{*} \in \mathcal{C}(\textbf{J},\mathbb{R}),$ we have
\begin{align*}
|(\textrm{\m{P}}\vartheta)(\textrm{\m{z}})-(\textrm{\m{P}}\vartheta^{*})(\textrm{\m{z}})|
\leq \quad &
\frac{1}{\Gamma(\alpha)}
\int_{0}^{\textrm{\m{z}}}{\psi}^{'}(\kappa)({\psi}(\textrm{\m{z}})-{\psi}(\kappa)^{\alpha-1}
\big|\mathscr{F}\big(\kappa, \vartheta(\kappa),\textrm{\m{I}}^{\alpha, \psi}_{0+}
     \mathscr{H}(\kappa, \tau, \textrm{\m{N}}^{\alpha, \psi}_{0+}\vartheta(\tau))\big)\\
       &-\mathscr{F}\big(\kappa, \vartheta^{*}(\kappa),\textrm{\m{I}}^{\alpha, \psi}_{0+}
    \mathscr{H}(\kappa, \tau, \textrm{\m{N}}^{\alpha, \psi}_{0+}\vartheta^{*}(\tau))\big)\big|d\kappa\\
\leq \quad & \frac{1}{\Gamma(\alpha)}
          \int_{0}^{\textrm{\m{z}}}{\psi}^{'}(\kappa)({\psi}(\textrm{\m{z}})-{\psi}(\kappa)^{\alpha-1}
             \textrm{\m{n}}(\kappa)
               \big[|\vartheta(\kappa)-\vartheta^{*}(\kappa)|\\
            &+|\textrm{\m{I}}^{\alpha, \psi}_{0+}
              \mathscr{H}(\kappa, \tau, \textrm{\m{N}}^{\alpha, \psi}_{0+}\vartheta(\tau))
               -\textrm{\m{I}}^{\alpha, \psi}_{0+}
            \mathscr{H}(\kappa, \tau, \textrm{\m{N}}^{\alpha, \psi}_{0+}\vartheta^{*}(\tau))|\big]d\kappa\\
\leq \quad & \frac{\mathscr{W}_{\textbf{1}}}{\Gamma(\alpha)}
              \int_{0}^{\textrm{\m{z}}}{\psi}^{'}(\kappa)({\psi}(\textrm{\m{z}})-{\psi}(\kappa)^{\alpha-1}
               \big[|\vartheta(\kappa)-\vartheta^{*}(\kappa)|\\
              &+\textrm{\m{I}}^{\alpha, \psi}_{0+}\textrm{\m{k}}(\kappa)
             |\textrm{\m{N}}^{\alpha, \psi}_{0+}\vartheta(\tau)
          -\textrm{\m{N}}^{\alpha, \psi}_{0+}\vartheta^{*}(\tau)|\big]d\kappa\\
\leq \quad & \frac{\mathscr{W}_{\textbf{1}}}{\Gamma(\alpha)}
             \int_{0}^{\textrm{\m{z}}}{\psi}^{'}(\kappa)({\psi}(\textrm{\m{z}})-{\psi}(\kappa)^{\alpha-1}
               \big[|\vartheta(\kappa)-\vartheta^{*}(\kappa)|\\
               &+\mathscr{W}_{\textbf{2}}\textrm{\m{I}}^{\alpha, \psi}_{0+}
             \textrm{\m{N}}^{\alpha, \psi}_{0+}|\vartheta(\tau)-\vartheta^{*}(\tau)|\big]d\kappa.\\
\end{align*}
Then in view of lemma (2.1), we obtain
\begin{align*}
|(\textrm{\m{P}}\vartheta)(\textrm{\m{z}})-(\textrm{\m{P}}\vartheta^{*})(\textrm{\m{z}})|
\leq \quad &\frac{\mathscr{W}_{\textbf{1}}}{\Gamma(\alpha)}
             \int_{0}^{\textrm{\m{z}}}{\psi}^{'}(\kappa)({\psi}(\textrm{\m{z}})-{\psi}(\kappa)^{\alpha-1}
            \big[|\vartheta(\kappa)-\vartheta^{*}(\kappa)|\\
            &+\mathscr{W}_{\textbf{2}}|\vartheta(\kappa)-\vartheta^{*}(\kappa)|\big]d\kappa,
\end{align*}
which gives
\[
\| (\textrm{\m{P}}\vartheta)(\textrm{\m{z}})-(\textrm{\m{P}}\vartheta^{*})(\textrm{\m{z}}) \|
\leq \frac{(\psi(b)-\psi(0))^{\alpha}}{\Gamma(\alpha+1)} \mathscr{W}_{\textbf{1}}
(1+\mathscr{W}_{\textbf{2}})\| \vartheta-\vartheta^{*}\|.
\]
By \eqref{eq3.1}, we have
\[
 \| (\textrm{\m{P}}\vartheta)(\textrm{\m{z}})-(\textrm{\m{P}}\vartheta^{*})(\textrm{\m{z}}) \|
\leq \| \vartheta-\vartheta^{*}\|.
\]
Hence $\textrm{\m{P}}$ is a contraction mapping and therefore, by Banach fixed point theorem $\textrm{\m{P}}$ has
unique fixed point in $\mathcal{C}(\textbf{J},\mathbb{R}),$ the fixed point is also a solution of \eqref{eq1.1}-\eqref{eq1.2}.$\qed$
\section{Estimate on the solution}
In this section, by using $\psi$-fractional Gronwall inequality, we prove estimate of
the solution of the problem \eqref{eq1.1}-\eqref{eq1.2}.\\

{\textbf{Theorem 4.1}}
Assume that the function $\mathscr{F}:\textbf{J}\times\mathbb{R}\times\mathbb{R}\rightarrow\mathbb{R}$ and $\mathscr{H}:\textbf{J}\times\textbf{J}\times\mathbb{R}\rightarrow\mathbb{R}$ satisfies the hypothesis ($\mathscr{Q1}$) and ($\mathscr{Q2}$)
respectively. If $\vartheta(\textrm{\m{z}}), \textrm{\m{z}}\in\textbf{J},$ is any solution of the problem
\eqref{eq1.1}-\eqref{eq1.2} then following holds:
\begin{align*}\label{eq4.1}
  |\vartheta(\textrm{\m{z}})|
\leq \quad &\big[|\vartheta_{0}|+\frac{(\psi(b)-\psi(0))^{\alpha}}{\Gamma(\alpha+1)}
(\mathscr{W}_{\textbf{3}}+\mathscr{W}_{\textbf{2}}\vartheta_{0})\big]\\
&\times E_{\alpha}\big(\mathscr{W}_{\textbf{1}}(1+\mathscr{W}_{\textbf{2}})(\psi(b)-\psi(0))^{\alpha}\big).\tag{4.1}
\end{align*}

{\textbf{Proof.}}
From lemma (2.6), we have
\begin{align*}
\vartheta({\textrm{\m{z}}})=\;\vartheta_{0}+\frac{1}{\Gamma(\alpha)}&
                            \int_{0}^{\textrm{\m{z}}}\psi^{'}(\kappa)(\psi(\textrm{\m{z}})-\psi(\kappa))^{\alpha-1}\\
                              &\mathscr{F}\big(\kappa, \vartheta(\kappa), \textrm{\m{I}}^{\alpha, \psi}_{0+}
                          \mathscr{H}(\kappa, \tau, \textrm{\m{N}}^{\alpha, \psi}_{0+}\vartheta(\tau))\big)d\kappa.
\end{align*}
Then, by hypothesis ($\mathscr{Q1}$), ($\mathscr{Q2}$) and for $\textrm{\m{z}}\in\textbf{J},$ we have
\begin{align*}
| \vartheta({\textrm{\m{z}}})| \leq \quad & | \vartheta_{0}|
                                 +\frac{1}{\Gamma(\alpha)}
                                  \int_{0}^{\textrm{\m{z}}}\psi^{'}(\kappa)(\psi(\textrm{\m{z}})-\psi(\kappa))^{\alpha-1}\\
                                       &|\mathscr{F}(\kappa, \vartheta(\kappa), \textrm{\m{I}}^{\alpha, \psi}_{0+}
                                          \mathscr{H}(\kappa, \tau, \textrm{\m{N}}^{\alpha, \psi}_{0+}\vartheta(\tau)))
                                       -\mathscr{F}(\kappa, 0, \textrm{\m{I}}^{\alpha, \psi}_{0+}\mathscr{H}(\kappa, \tau, 0))|d\kappa\\
                                      +&\frac{1}{\Gamma(\alpha)}
                                 \int_{0}^{\textrm{\m{z}}}\psi^{'}(\kappa)(\psi(\textrm{\m{z}})-\psi(\kappa))^{\alpha-1}
                               | \mathscr{F}(\kappa, 0, \textrm{\m{I}}^{\alpha, \psi}_{0+} \mathscr{H}(\kappa, \tau, 0))|d\kappa \\
\leq \quad & | \vartheta_{0}| +\frac{\mathscr{W}_{\textbf{1}}}{\Gamma(\alpha)}
                    \int_{0}^{\textrm{\m{z}}}\psi^{'}(\kappa)(\psi(\textrm{\m{z}})-\psi(\kappa))^{\alpha-1}\\
                        &[|\vartheta(\kappa)|+\textrm{\m{I}}^{\alpha, \psi}_{0+}
                          |\mathscr{H}(\kappa, \tau, \textrm{\m{N}}^{\alpha, \psi}_{0+}\vartheta(\tau))-\mathscr{H}(\kappa, \tau, 0)|]d\kappa\\
                     +&\frac{(\psi(b)-\psi(0))^{\alpha}}{\Gamma(\alpha+1)}\mathscr{W}_{\textbf{3}} \\
\leq \quad & | \vartheta_{0}| +\frac{\mathscr{W}_{\textbf{1}}}{\Gamma(\alpha)}
                    \int_{0}^{\textrm{\m{z}}}\psi^{'}(\kappa)(\psi(\textrm{\m{z}})-\psi(\kappa))^{\alpha-1}\\
                        &[|\vartheta(\kappa)|+\mathscr{W}_{\textbf{2}}
                        \textrm{\m{I}}^{\alpha, \psi}_{0+}\textrm{\m{N}}^{\alpha, \psi}_{0+}|\vartheta(\tau)|]d\kappa
                     +\frac{(\psi(b)-\psi(0))^{\alpha}}{\Gamma(\alpha+1)}\mathscr{W}_{\textbf{3}}.
\end{align*}
Then, by lemma (2.1), we obtain
\begin{align*}\label{eq4.2}
| \vartheta({\textrm{\m{z}}})| \leq \quad & | \vartheta_{0}| +\frac{\mathscr{W}_{\textbf{1}}}{\Gamma(\alpha)}
                    \int_{0}^{\textrm{\m{z}}}\psi^{'}(\kappa)(\psi(\textrm{\m{z}})-\psi(\kappa))^{\alpha-1}\\
                        &[|\vartheta(\kappa)|+\mathscr{W}_{\textbf{2}}
                        |\vartheta(\kappa)-\vartheta_{0}|]d\kappa
                     +\frac{(\psi(b)-\psi(0))^{\alpha}}{\Gamma(\alpha+1)}\mathscr{W}_{\textbf{3}}\\
= \quad & | \vartheta_{0}| +\mathscr{W}_{\textbf{1}}(1+\mathscr{W}_{\textbf{2}})
                        \int_{0}^{\textrm{\m{z}}}\psi^{'}(\kappa)(\psi(\textrm{\m{z}})-\psi(\kappa))^{\alpha-1}
                        \vartheta(\kappa)d\kappa\\
                     &+\frac{(\psi(b)-\psi(0))^{\alpha}}{\Gamma(\alpha+1)}
                     \big(\mathscr{W}_{\textbf{3}}+\mathscr{W}_{\textbf{2}}\vartheta_{0}\big).\tag{4.2}
\end{align*}
Applying lemma (2.5) to \eqref{eq4.2}, we get
\begin{align*}
 | \vartheta({\textrm{\m{z}}})| \leq \quad &
   | \vartheta_{0}|+\frac{(\psi(b)-\psi(0))^{\alpha}}{\Gamma(\alpha+1)}
     \big(\mathscr{W}_{\textbf{3}}+\mathscr{W}_{\textbf{2}}\vartheta_{0}\big)\\
      &+\int_{0}^{\textrm{\m{z}}}\sum_{m=1}^{\infty}
       \frac{(\mathscr{W}_{\textbf{1}}(1+\mathscr{W}_{\textbf{2}}))^{m}}{\Gamma(\alpha m)}
      \psi^{'}(\kappa)(\psi(\textrm{\m{z}})-\psi(\kappa))^{\alpha m-1}\\
    &\Big\{\vartheta_{0}+\Big(\frac{(\psi(b)-\psi(0))^{\alpha}}{\Gamma(\alpha+1)}\Big)
  \big(\mathscr{W}_{\textbf{3}}+\mathscr{W}_{\textbf{2}}\vartheta_{0}\big)\Big\}d\kappa.\\
 \end{align*}
This gives
\begin{align*}
| \vartheta({\textrm{\m{z}}})| \leq \quad & \big[| \vartheta_{0}|+\frac{(\psi(b)-\psi(0))^{\alpha}}{\Gamma(\alpha+1)}
     \big(\mathscr{W}_{\textbf{3}}+\mathscr{W}_{\textbf{2}}\vartheta_{0}\big)\big]\\
     &\Big[1+\int_{0}^{\textrm{\m{z}}}\sum_{m=1}^{\infty}
       \frac{(\mathscr{W}_{\textbf{1}}(1+\mathscr{W}_{\textbf{2}}))^{m}}{\Gamma(\alpha m)}
      \psi^{'}(\kappa)(\psi(\textrm{\m{z}})-\psi(\kappa))^{\alpha m-1}\Big]d\kappa\\
= \quad &   \big[| \vartheta_{0}|+\frac{(\psi(b)-\psi(0))^{\alpha}}{\Gamma(\alpha+1)}
     \big(\mathscr{W}_{\textbf{3}}+\mathscr{W}_{\textbf{2}}\vartheta_{0}\big)\big]\\
      &\Big[1+\sum_{m=1}^{\infty}\frac{\big[\mathscr{W}_{\textbf{1}}(1+\mathscr{W}_{\textbf{2}})
        (\psi(b)-\psi(0))^{\alpha}\big]^{m}}{\Gamma(\alpha m+1)}\Big]\\
= \quad &  \big[ | \vartheta_{0}|+\frac{(\psi(b)-\psi(0))^{\alpha}}{\Gamma(\alpha+1)}
     \big(\mathscr{W}_{\textbf{3}}+\mathscr{W}_{\textbf{2}}\vartheta_{0}\big)\big]\\
    & E_{\alpha}(\mathscr{W}_{\textbf{1}}(1+\mathscr{W}_{\textbf{2}})(\psi(b)-\psi(0))^{\alpha}).\qed
\end{align*}
\section{Stability Analysis}
In this section, we investigate the Ulam stability results for the problem \eqref{eq1.1}-\eqref{eq1.2}.\\
{\textbf{Definition 5.1}}.\cite{Abb}.
The problem \eqref{eq1.1}-\eqref{eq1.2} is said to be Ulam-Hyers(UH) stable if there exists a real number $\gamma>0$
such that for every $\varepsilon>0$ and for each solution $\omega\in\mathcal{C}(\textbf{J},\mathbb{R})$ of the inequality
\[\label{eq5.1}
\Big|\textrm{\m{N}}^{\alpha, \psi}_{0+}\omega(\textrm{\m{z}})
-\mathscr{F}
  \big(\textrm{\m{z}}, \omega(\textrm{\m{z}}),
    \textrm{\m{I}}^{\alpha, \psi}_{0+}
    \mathscr{H}
       (\textrm{\m{z}}, \tau, \textrm{\m{N}}^{\alpha, \psi}_{0+}\omega(\tau))\big)\Big|\leq\varepsilon,
\tag{5.1}
\]
there exists a solution $\vartheta\in\mathcal{C}(\textbf{J},\mathbb{R})$ of the problem \eqref{eq1.1}-\eqref{eq1.2} with
\[
|\omega(\textrm{\m{z}})-\vartheta(\textrm{\m{z}})|\leq\gamma\varepsilon,\quad \textrm{\m{z}}\in [0, b].
\]
{\textbf{Definition 5.2}}.\cite{Abb}
The problem \eqref{eq1.1}-\eqref{eq1.2} is said to be Generalized Ulam-Hyers(GUH) stable if there exists a continuous function
$\rho:\mathbb{R}^{+}\rightarrow\mathbb{R}^{+}$ with $\rho(0)=0$ such that for every $\varepsilon>0$ and for each solution $\omega\in\mathcal{C}(\textbf{J},\mathbb{R})$ of the inequality \eqref{eq5.1},
there exist a solution $\vartheta\in\mathcal{C}(\textbf{J},\mathbb{R})$ of the problem \eqref{eq1.1}-\eqref{eq1.2} with
\[
|\omega(\textrm{\m{z}})-\vartheta(\textrm{\m{z}})|\leq\rho(\varepsilon),\quad \textrm{\m{z}}\in [0, b].
\]
{\textbf{Remark 5.1}}
A function $\omega\in\mathcal{C}(\textbf{J},\mathbb{R})$ is a solution of the inequality \eqref{eq5.1} if and only if there exists a
function $\textrm{\B{h}}\in\mathcal{C}(\textbf{J},\mathbb{R})$ (where $\textrm{\B{h}}$ depends on $\omega$) such that
\begin{align*}
&(1)|\textrm{\B{h}}(\textrm{\m{z}})|<\varepsilon\\
&(2)\textrm{\m{N}}^{\alpha, \psi}_{0+}\omega(\textrm{\m{z}})
=\mathscr{F}
  \big(\textrm{\m{z}}, \omega(\textrm{\m{z}}),
    \textrm{\m{I}}^{\alpha, \psi}_{0+}
    \mathscr{H} (\textrm{\m{z}}, \tau, \textrm{\m{N}}^{\alpha, \psi}_{0+}\omega(\tau))\big)+\textrm{\B{h}}(\textrm{\m{z}}),\;\textrm{\m{z}}\in[0,b].
\end{align*}
{\textbf{Definition 5.3}.\cite{Abb}
The problem \eqref{eq1.1}-\eqref{eq1.2} is said to be Ulam-Hyers-Rassias(UHR) stable with respect to
$\rho \in \mathcal{C}(\textbf{J},\mathbb{R}),$ if there exists a real number $B>0$ such that
for every $\varepsilon>0$ and for each solution $\omega\in\mathcal{C}(\textbf{J},,\mathbb{R})$ of the inequality
\[\label{eq5.2}
\Big|\textrm{\m{N}}^{\alpha, \psi}_{0+}\omega(\textrm{\m{z}})
-\mathscr{F}
  \big(\textrm{\m{z}}, \omega(\textrm{\m{z}}),
    \textrm{\m{I}}^{\alpha, \psi}_{0+}
    \mathscr{H}
       (\textrm{\m{z}}, \tau, \textrm{\m{N}}^{\alpha, \psi}_{0+}\omega(\tau))\big)\Big|
         \leq \varepsilon\rho(\textrm{\m{z}}),
\tag{5.2}
\]
there exist a solution $\vartheta\in\mathcal{C}(\textbf{J},\mathbb{R})$ of the problem \eqref{eq1.1}-\eqref{eq1.2} with
\[
|\omega(\textrm{\m{z}})-\vartheta(\textrm{\m{z}})|\leq B\rho(\textrm{\m{z}})\varepsilon,\quad \textrm{\m{z}}\in [0, b].
\]
{\textbf{Definition 5.4}.\cite{Abb}
The problem \eqref{eq1.1}-\eqref{eq1.2} is said to be Generalized Ulam-Hyers-Rassias(GUHR) stable with respect to
$\rho \in \mathcal{C}(\textbf{J},\mathbb{R}),$ if there exists a real number $B>0$ such that
for each solution $\omega\in\mathcal{C}(\textbf{J},\mathbb{R})$ of the inequality
\[\label{eq5.3}
\Big|\textrm{\m{N}}^{\alpha, \psi}_{0+}\omega(\textrm{\m{z}})
-\mathscr{F}
  \big(\textrm{\m{z}}, \omega(\textrm{\m{z}}),
    \textrm{\m{I}}^{\alpha, \psi}_{0+}
    \mathscr{H}
       \big(\textrm{\m{z}}, \tau, \textrm{\m{N}}^{\alpha, \psi}_{0+}\omega(\tau)\big)\big)\Big|
         \leq \rho(\textrm{\m{z}}),
\tag{5.3}
\]
there exist a solution $\vartheta\in\mathcal{C}(\textbf{J},\mathbb{R})$ of the problem \eqref{eq1.1}-\eqref{eq1.2} with
\[
|\omega(\textrm{\m{z}})-\vartheta(\textrm{\m{z}})|\leq B\rho(\textrm{\m{z}}),\quad \textrm{\m{z}}\in [0, b].
\]
{\textbf{Remark 5.2}}.
A function $\omega\in\mathcal{C}(\textbf{J},\mathbb{R})$ is a solution of the inequality \eqref{eq6.2} if and only if there exists
function $\textrm{\B{h}}, \rho\in\mathcal{C}(\textbf{J},\mathbb{R})$ (where $\textrm{\B{h}}$ depends on $\omega$) such that
\begin{align*}
&(1)|\textrm{\B{h}}(\textrm{\m{z}})|<\varepsilon\rho(\textrm{\m{z}})\\
&(2)\textrm{\m{N}}^{\alpha, \psi}_{0+}\omega(\textrm{\m{z}})
=\mathscr{F}
  \big(\textrm{\m{z}}, \omega(\textrm{\m{z}}),
    \textrm{\m{I}}^{\alpha, \psi}_{0+}
    \mathscr{H}  \big(\textrm{\m{z}}, \tau, \textrm{\m{N}}^{\alpha, \psi}_{0+}\omega(\tau)\big)\big)+\textrm{\B{h}}(\textrm{\m{z}}),\;\textrm{\m{z}}\in[0,b].
\end{align*}
{\textbf{Theorem 5.1}} Assume that the hypothesis of Theorem (3.1) hold. Let the inequality \eqref{eq5.1} is satisfied,
then the problem \eqref{eq1.1}-\eqref{eq1.1} is UH stable.\\
{\textbf{Proof.}}
Let $\omega \in \mathcal{C}(\textbf{J},\mathbb{R})$ be a solution of inequality \eqref{eq5.1} and $\vartheta$ be a solution of
\[\label{eq5.4}
\textrm{\m{N}}^{\alpha, \psi}_{0+}\vartheta(\textrm{\m{z}})
=\mathscr{F}
  \Big(\textrm{\m{z}}, \vartheta(\textrm{\m{z}}),
    \textrm{\m{I}}^{\alpha, \psi}_{0+}
    \mathscr{H}
       \big(\textrm{\m{z}}, \tau, \textrm{\m{N}}^{\alpha, \psi}_{0+}\vartheta(\tau)\big)\Big),
       \;\mathrm{for\; all}\; \textrm{\m{z}}\in \textbf{J}, 0<\alpha<1,
\tag{5.4}
\]
\[\label{eq5.5}
\vartheta(0)=\omega(0)=\vartheta_{0}=\omega_{0}
\tag{5.5}
\]
Then we have
\begin{align*}
\vartheta({\textrm{\m{z}}})= \vartheta_{0}+\frac{1}{\Gamma(\alpha)}&
                            \int_{0}^{\textrm{\m{z}}}\psi^{'}(\kappa)(\psi(\textrm{\m{z}})-\psi(\kappa))^{\alpha-1}\\
                              &\mathscr{F}\big(\kappa, \vartheta(\kappa), \textrm{\m{I}}^{\alpha, \psi}_{0+}
                          \mathscr{H}(\kappa, \tau, \textrm{\m{N}}^{\alpha, \psi}_{0+}\vartheta(\tau))\big)d\kappa.
\end{align*}
From \eqref{eq5.5}, $\vartheta_{0}= \omega_{0}$ and therefore, we get
\begin{align*}\label{eq5.6}
\vartheta({\textrm{\m{z}}})= \omega_{0}+\frac{1}{\Gamma(\alpha)}&
                            \int_{0}^{\textrm{\m{z}}}\psi^{'}(\kappa)(\psi(\textrm{\m{z}})-\psi(\kappa))^{\alpha-1}\\
                              &\mathscr{F}\big(\kappa, \vartheta(\kappa), \textrm{\m{I}}^{\alpha, \psi}_{0+}
                          \mathscr{H}(\kappa, \tau, \textrm{\m{N}}^{\alpha, \psi}_{0+}\vartheta(\tau))\big)d\kappa.
      \tag{5.6}
\end{align*}
Since $\omega(\textrm{\m{z}})$ is a solution of the inequality \eqref{eq5.1}, then by Remark (5.1), there exist a function
$\textrm{\B{h}}(\textrm{\m{z}})$ such that
\[\label{eq5.7}
|\textrm{\B{h}}(\textrm{\m{z}})|\leq\varepsilon
\tag{5.7}
\]
and
\[\label{eq5.8}
\textrm{\m{N}}^{\alpha, \psi}_{0+}\omega(\textrm{\m{z}})
=\mathscr{F}
  \big(\textrm{\m{z}}, \omega(\textrm{\m{z}}),
    \textrm{\m{I}}^{\alpha, \psi}_{0+}
    \mathscr{H}  \big(\textrm{\m{z}}, \tau, \textrm{\m{N}}^{\alpha, \psi}_{0+}\omega(\tau)\big)\big)+\textrm{\B{h}}(\textrm{\m{z}}),\;\textrm{\m{z}}\in [0,b].
\tag{5.8}\]
From \eqref{eq5.7} and \eqref{eq5.8}, we obtain
\begin{align*}\label{eq5.9}
\Big|\omega(\textrm{\m{z}})-\omega_{0}-\frac{1}{\Gamma(\alpha)}&
           \int_{0}^{\textrm{\m{z}}}\psi^{'}(\kappa)(\psi(\textrm{\m{z}})-\psi(\kappa))^{\alpha-1}\\
          & \mathscr{F}\big(\kappa, \omega(\kappa), \textrm{\m{I}}^{\alpha, \psi}_{0+}
           \mathscr{H}(\kappa, \tau, \textrm{\m{N}}^{\alpha, \psi}_{0+}\omega(\tau))\big)d\kappa\Big|
     \leq \frac{(\psi(b)-\psi(0))^{\alpha}}{\Gamma(\alpha+1)} \varepsilon.\tag{5.9}
\end{align*}
Using hypothesis ($\mathscr{Q1}$), ($\mathscr{Q2}$) and from \eqref{eq5.6} and \eqref{eq5.9}, we obtain
\begin{align*}\label{eq5.10}
|\omega(\textrm{\m{z}})-\vartheta(\textrm{\m{z}})|
 \leq \quad &\Big|\omega(\textrm{\m{z}})-\omega_{0}
   -\frac{1}{\Gamma(\alpha)}\int_{0}^{\textrm{\m{z}}}\psi^{'}(\kappa)(\psi(\textrm{\m{z}})-\psi(\kappa))^{\alpha-1}\\
    &\mathscr{F}\big(\kappa, \omega(\kappa), \textrm{\m{I}}^{\alpha, \psi}_{0+}
           \mathscr{H}(\kappa, \tau, \textrm{\m{N}}^{\alpha, \psi}_{0+}\omega(\tau))\big)d\kappa\Big|\\
        +\quad&\frac{1}{\Gamma(\alpha)}\int_{0}^{\textrm{\m{z}}}\psi^{'}(\kappa)(\psi(\textrm{\m{z}})-\psi(\kappa))^{\alpha-1}
    |\mathscr{F}\big(\kappa, \omega(\kappa), \textrm{\m{I}}^{\alpha, \psi}_{0+}
           \mathscr{H}(\kappa, \tau, \textrm{\m{N}}^{\alpha, \psi}_{0+}\omega(\tau))\big)\\
    -\quad&\mathscr{F}\big(\kappa, \vartheta(\kappa), \textrm{\m{I}}^{\alpha, \psi}_{0+}
           \mathscr{H}(\kappa, \tau, \textrm{\m{N}}^{\alpha, \psi}_{0+}\vartheta(\tau))\big)|d\kappa\\
 \leq \quad& \frac{(\psi(b)-\psi(0))^{\alpha}}{\Gamma(\alpha+1)} \varepsilon
 +\frac{\mathscr{W}_{\textbf{1}}  (1+\mathscr{W}_{\textbf{2}})}{\Gamma(\alpha)}
 \int_{0}^{\textrm{\m{z}}}\psi^{'}(\kappa)(\psi(\textrm{\m{z}})-\psi(\kappa))^{\alpha-1}\\
 \quad&|\omega(\kappa)-\vartheta(\kappa)|d\kappa.\tag{5.10}
 \end{align*}
 Applying lemma (2.5) to \eqref{eq5.10}, we get
\begin{align*}
|\omega(\textrm{\m{z}})-\vartheta(\textrm{\m{z}})|\\
\leq \quad&  \frac{(\psi(b)-\psi(0))^{\alpha}}{\Gamma(\alpha+1)} \varepsilon
     \Big[1+\int_{0}^{\textrm{\m{z}}}\sum_{m=1}^{\infty}
        \frac{[\mathscr{W}_{\textbf{1}}  (1+\mathscr{W}_{\textbf{2}})]^{m}}{\Gamma(\alpha m)}
      \psi^{'}(\kappa)(\psi(\textrm{\m{z}})-\psi(\kappa))^{\alpha m-1}d\kappa\Big]\\
= \quad& \frac{(\psi(b)-\psi(0))^{\alpha}}{\Gamma(\alpha+1)} \varepsilon
 \Big[1+\sum_{m=1}^{\infty}\frac{[\mathscr{W}_{\textbf{1}}  (1+\mathscr{W}_{\textbf{2}})]^{m}}
     {\Gamma(\alpha m)}\int_{0}^{\textrm{\m{z}}}\psi^{'}(\kappa)(\psi(\textrm{\m{z}})-\psi(\kappa))^{\alpha m-1}d\kappa\Big]\\
\leq \quad &\frac{(\psi(b)-\psi(0))^{\alpha}}{\Gamma(\alpha+1)} \varepsilon
     \Big[1+\sum_{m=1}^{\infty}\frac{[\mathscr{W}_{\textbf{1}}  (1+\mathscr{W}_{\textbf{2}})]^{m}}
        {\Gamma(\alpha m+1)}(\psi(b)-\psi(0))^{\alpha m}\Big]\\
= \quad & \frac{(\psi(b)-\psi(0))^{\alpha}\varepsilon}{\Gamma(\alpha+1)}
    E_{\alpha}(\mathscr{W}_{\textbf{1}}  (1+\mathscr{W}_{\textbf{2}})(\psi(b)-\psi(0))^{\alpha})).
\end{align*}
Put
\[
\gamma=\frac{(\psi(b)-\psi(0))^{\alpha}}{\Gamma(\alpha+1)}
    E_{\alpha}(\mathscr{W}_{\textbf{1}}  (1+\mathscr{W}_{\textbf{2}})(\psi(b)-\psi(0))^{\alpha})).
\]
It follows that
\[
|\omega(\textrm{\m{z}})-\vartheta(\textrm{\m{z}})|\leq \gamma\varepsilon.
\]
Hence, the problem \eqref{eq1.1}-\eqref{eq1.2} is UH stable.$\qed$\\

{\textbf{Theorem 5.2}} Assume that the hypothesis of Theorem (5.1) hold and there exists a function $\rho\in\mathcal{C}(\mathbb{R}^{+},\mathbb{R}^{+})$ such that $\rho(0)=0.$ Then the problem
\eqref{eq1.1}-\eqref{eq1.2} is GUH stable.\\
{\textbf{Proof.}}
Put $\rho(\varepsilon)=\gamma \varepsilon$ with $\rho(0)=0$ and using same argument as in Theorem (5.1), we obtain
\[
|\omega(\textrm{\m{z}})-\vartheta(\textrm{\m{z}})|\leq \rho(\varepsilon). \qed
\]

{\textbf{Lemma 5.1}} Consider the following hypothesis:\\
($\mathscr{Q3}$)There exist an increasing function $\rho(\textrm{\m{z}}) \in \mathcal{C}(\textbf{J},\mathbb{R})$ and $\gamma>0$ such that
\[
\textrm{\m{I}}^{\alpha, \psi}_{0+}\rho(\textrm{\m{z}})\leq\gamma\rho(\textrm{\m{z}}),\quad\textrm{\m{z}}\in[0,b].
\]
Let $\omega(\textrm{\m{z}}) \in \mathcal{C}(\textbf{J},\mathbb{R})$ be a solution of the inequality \eqref{eq5.2}.
Then following hold
\begin{align*}
\big|\omega(\textrm{\m{z}})-\omega_{0}-&\frac{1}{\Gamma(\alpha)}
                                \int_{0}^{\textrm{\m{z}}}\psi^{'}(\kappa)(\psi(\textrm{\m{z}})-\psi(\kappa))^{\alpha-1}\\
                                    &\mathscr{F}(\kappa, \vartheta(\kappa), \textrm{\m{I}}^{\alpha, \psi}_{0+}
                               \mathscr{H}(\kappa, \tau, \textrm{\m{N}}^{\alpha, \psi}_{0+}\vartheta(\tau))))d\kappa\big|
\leq \varepsilon\gamma\rho(\textrm{\m{z}}).
\end{align*}
{\textbf{Proof}} Let $\varepsilon>0$ and $\omega$ be a solution of the inequality \eqref{eq5.2}, then by Remark (5.2),
there exists a function $\textrm{\B{h}}, \rho \in \mathcal{C}(\textbf{J},\mathbb{R})$ such that
\[
|\textrm{\B{h}}(\textrm{\m{z}})|\leq\varepsilon\rho(\textrm{\m{z}})
\]
and
\[
\textrm{\m{N}}^{\alpha, \psi}_{0+}\omega(\textrm{\m{z}})
=\mathscr{F}
  \big(\textrm{\m{z}}, \omega(\textrm{\m{z}}),
    \textrm{\m{I}}^{\alpha, \psi}_{0+}
    \mathscr{H}  \big(\textrm{\m{z}}, \tau, \textrm{\m{N}}^{\alpha, \psi}_{0+}\omega(\tau)\big)\big)+\textrm{\B{h}}(\textrm{\m{z}}),\;\textrm{\m{z}}\in [0,b].
\]
It follows that
\begin{align*}
\omega(\textrm{\m{z}})&=\omega_{0}+\frac{1}{\Gamma(\alpha)}
                                \int_{0}^{\textrm{\m{z}}}\psi^{'}(\kappa)(\psi(\textrm{\m{z}})-\psi(\kappa))^{\alpha-1}
                                    \mathscr{F}(\kappa, \omega(\kappa), \textrm{\m{I}}^{\alpha, \psi}_{0+}
                               \mathscr{H}(\kappa, \tau, \textrm{\m{N}}^{\alpha, \psi}_{0+}\omega(\tau)))d\kappa\\
                            &+\frac{1}{\Gamma(\alpha)}
                          \int_{0}^{\textrm{\m{z}}}\psi^{'}(\kappa)(\psi(\textrm{\m{z}})-\psi(\kappa))^{\alpha-1}
                             \textrm{\B{h}}(\kappa)d\kappa.
\end{align*}
By using hypothesis ($\mathscr{Q3}$), we get
\begin{align*}
\Big|\omega(\textrm{\m{z}})-&\omega_{0}-\frac{1}{\Gamma(\alpha)}
                                \int_{0}^{\textrm{\m{z}}}\psi^{'}(\kappa)(\psi(\textrm{\m{z}})-\psi(\kappa))^{\alpha-1}
                                    \mathscr{F}(\kappa, \omega(\kappa), \textrm{\m{I}}^{\alpha, \psi}_{0+}
                               \mathscr{H}(\kappa, \tau, \textrm{\m{N}}^{\alpha, \psi}_{0+}\omega(\tau)))d\kappa\Big|\\
 \leq \quad&  \frac{1}{\Gamma(\alpha)}
             \int_{0}^{\textrm{\m{z}}}\psi^{'}(\kappa)(\psi(\textrm{\m{z}})-\psi(\kappa))^{\alpha-1}
         |\textrm{\B{h}}(\kappa)|d\kappa\\
 \leq \quad& \frac{\varepsilon}{\Gamma(\alpha)}
            \int_{0}^{\textrm{\m{z}}}\psi^{'}(\kappa)(\psi(\textrm{\m{z}})-\psi(\kappa))^{\alpha-1}
         \rho(\kappa)d\kappa\\
 \leq \quad& \varepsilon \gamma \rho(\textrm{\m{z}}). \qed\\
\end{align*}

{\textbf{Theorem 5.3}} Suppose that hypothesis ($\mathscr{Q1}$), ($\mathscr{Q2}$) and ($\mathscr{Q3}$) holds. Let $\omega\in\mathcal{C}(\textbf{J},\mathbb{R})$
be a solution of the inequality \eqref{eq5.2} and $\mathscr{W}_{\textbf{1}}(1+\mathscr{W}_{\textbf{2}})\gamma\neq1,$
then the problem \eqref{eq1.1}-\eqref{eq1.2} is UHR stable.\\
{\textbf{Proof}} Let $\varepsilon>0$ and $\vartheta (\textrm{\m{z}})$ be a solution of the problem \eqref{eq1.1}-\eqref{eq1.2}.
Then we have
\begin{align*}
\vartheta(\textrm{\m{z}})=\omega_{0}+&\frac{1}{\Gamma(\alpha)}
                                \int_{0}^{\textrm{\m{z}}}\psi^{'}(\kappa)(\psi(\textrm{\m{z}})-\psi(\kappa))^{\alpha-1}\\
                                    &\mathscr{F}(\kappa, \vartheta(\kappa), \textrm{\m{I}}^{\alpha, \psi}_{0+}
                               \mathscr{H}(\kappa, \tau, \textrm{\m{N}}^{\alpha, \psi}_{0+}\vartheta(\tau)))d\kappa.
\end{align*}
By hypothesis ($\mathscr{Q1}$), ($\mathscr{Q2}$) and Lemma (5.1), we obtain
\begin{align*}\label{eq5.11}
|\omega(\textrm{\m{z}})-\vartheta(\textrm{\m{z}})|
\leq \quad& \Big|\omega(\textrm{\m{z}})-\omega_{0}-\frac{1}{\Gamma(\alpha)}
      \int_{0}^{\textrm{\m{z}}} \psi^{'}(\kappa)(\psi(\textrm{\m{z}})-\psi(\kappa))^{\alpha-1}\\
        &\quad  \mathscr{F}(\kappa, \omega(\kappa), \textrm{\m{I}}^{\alpha, \psi}_{0+}
              \mathscr{H}(\kappa, \tau, \textrm{\m{N}}^{\alpha, \psi}_{0+}\omega(\tau)))d\kappa\Big| \\
  + \quad&\frac{1}{\Gamma(\alpha)}\int_{0}^{\textrm{\m{z}}} \psi^{'}(\kappa)(\psi(\textrm{\m{z}})-\psi(\kappa))^{\alpha-1}
    |\mathscr{F}(\kappa, \omega(\kappa), \textrm{\m{I}}^{\alpha, \psi}_{0+}
              \mathscr{H}(\kappa, \tau, \textrm{\m{N}}^{\alpha, \psi}_{0+}\omega(\tau)))\\
 &\quad- \mathscr{F}(\kappa, \vartheta(\kappa), \textrm{\m{I}}^{\alpha, \psi}_{0+}
              \mathscr{H}(\kappa, \tau, \textrm{\m{N}}^{\alpha, \psi}_{0+}\vartheta(\tau)))|d\kappa\\
\leq \quad& \varepsilon\gamma\rho(\textrm{\m{z}})
         +\frac{\mathscr{W}_{\textbf{1}}(1+\mathscr{W}_{\textbf{2}})}{\Gamma(\alpha)}
    \int_{0}^{\textrm{\m{z}}} \psi^{'}(\kappa)(\psi(\textrm{\m{z}})-\psi(\kappa))^{\alpha-1}
       |\omega(\kappa)-\vartheta(\kappa)|d\kappa.\tag{5.11}
\end{align*}
Applying Lemma (2.5) to \eqref{eq5.11} and using hypothesis ($\mathscr{Q3}$), we obtain
\begin{align*}
|\omega(\textrm{\m{z}})-\vartheta(\textrm{\m{z}})|
\leq\quad& \varepsilon\gamma\rho(\textrm{\m{z}})+\varepsilon\gamma
              \int_{0}^{\textrm{\m{z}}} \sum_{m=1}^{\infty}
           \frac{\big[\mathscr{W}_{\textbf{1}}(1+\mathscr{W}_{\textbf{2}})\big]^{m}}{\Gamma (\alpha m)}
         \psi^{'}(\kappa)(\psi(\textrm{\m{z}})-\psi(\kappa))^{\alpha m-1}\rho(\kappa)d\kappa\\
 =\quad& \varepsilon\gamma\rho(\textrm{\m{z}})+\varepsilon\gamma
           \Big[\int_{0}^{\textrm{\m{z}}}
              \frac{\mathscr{W}_{\textbf{1}}(1+\mathscr{W}_{\textbf{2}})}{\Gamma(\alpha)}
              \psi^{'}(\kappa)(\psi(\textrm{\m{z}})-\psi(\kappa))^{\alpha-1}\rho(\kappa)d\kappa\\
 +\quad & \int_{0}^{\textrm{\m{z}}}
              \frac{[\mathscr{W}_{\textbf{1}}(1+\mathscr{W}_{\textbf{2}})]^{2}}{\Gamma(2\alpha)}
              \psi^{'}(\kappa)(\psi(\textrm{\m{z}})-\psi(\kappa))^{2\alpha-1}\rho(\kappa)d\kappa+....\Big]\\
=\quad & \varepsilon\gamma\rho(\textrm{\m{z}})+\varepsilon\gamma
         \big[\mathscr{W}_{\textbf{1}}(1+\mathscr{W}_{\textbf{2}})
              \textrm{\m{I}}^{\alpha, \psi}_{0+}\rho(\textrm{\m{z}})
              +(\mathscr{W}_{\textbf{1}}(1+\mathscr{W}_{\textbf{2}}))^{2}
              \textrm{\m{I}}^{2\alpha, \psi}_{0+}\rho(\textrm{\m{z}})+....\big]\\
\leq \quad & \varepsilon\gamma\rho(\textrm{\m{z}})+\varepsilon\gamma
       \big[ \mathscr{W}_{\textbf{1}}(1+\mathscr{W}_{\textbf{2}})\gamma\rho(\textrm{\m{z}})
       +(\mathscr{W}_{\textbf{1}}(1+\mathscr{W}_{\textbf{2}})\gamma)^{2}\rho(\textrm{\m{z}})+....\big]\\
=\quad &\varepsilon\gamma\rho(\textrm{\m{z}})\sum_{m=0}^{\infty}
       (\mathscr{W}_{\textbf{1}}(1+\mathscr{W}_{\textbf{2}})\gamma)^{m}\\
= \quad &\frac{\varepsilon\gamma\rho(\textrm{\m{z}})}
{1-(\mathscr{W}_{\textbf{1}}(1+\mathscr{W}_{\textbf{2}}))\gamma}.
\end{align*}
Put
\[
B=\frac{\gamma}
{1-(\mathscr{W}_{\textbf{1}}(1+\mathscr{W}_{\textbf{2}}))\gamma}.
\]
It follows that
\[
|\omega(\textrm{\m{z}})-\vartheta(\textrm{\m{z}})|\leq B\varepsilon\rho(\textrm{\m{z}}), \quad \textrm{\m{z}}\in[0, b].
\]
Hence, the problem \eqref{eq1.1}-\eqref{eq1.2} is UHR stable. $\qed$\\

 {\textbf{Theorem 5.4}} Suppose that hypothesis ($\mathscr{Q1}$), ($\mathscr{Q2}$) and ($\mathscr{Q3}$) holds. Let $\omega\in\mathcal{C}(\textbf{J},\mathbb{R})$
be a solution of the inequality \eqref{eq5.3} and $\mathscr{W}_{\textbf{1}}(1+\mathscr{W}_{\textbf{2}})\gamma\neq1,$
then the problem \eqref{eq1.1}-\eqref{eq1.2} is GUHR stable.\\
{\textbf{Proof.}}
Put $\varepsilon=1$ and using same argument as in Theorem (5.3), we obtain
\[
|\omega(\textrm{\m{z}})-\vartheta(\textrm{\m{z}})|\leq B\rho(\textrm{\m{z}}).\qed
\]
\section{Concluding Remark}
In this study, applying the method of Banach fixed point theorem, we obtain existence and uniqueness of the fractional differential equation with $\psi$-Caputo derivative. Furthermore, estimate and some stability results for corresponding fractional differential equation are studied by using $\psi$-Gronwall inequality. For $\psi(\textrm{\m{z}})=\textrm{\m{z}},$ $ln(\textrm{\m{z}}),$ $\textrm{\m{z}}^{\sigma},$ $\psi$-Caputo derivative reduces to the Caputo \cite{Kil}, Caputo-Hadamard \cite{Gam}, Caputo–Erd\'{e}lyi–Kober \cite{Luc} fractional derivative respectively.

\end{document}